# Sample path properties of the local time of multifractional Brownian motion

BRAHIM BOUFOUSSI[1,*], MARCO DOZZI[2] and RABY GUERBAZ[1,**]

[1]*Department of Mathematics, Faculty of Sciences Semlalia, Cadi Ayyad University, 2390 Marrakesh, Morocco.* E-mail: [*]*boufoussi@ucam.ac.ma;* [**]*r.guerbaz@ucam.ac.ma*
[2]*IECN UMR 7502, Nancy Université, 54506 Vandoeuvre-Lès-Nancy, France.*
E-mail: *dozzi@iecn.u-nancy.fr*

We establish estimates for the local and uniform moduli of continuity of the local time of multifractional Brownian motion, $B^H = (B^{H(t)}(t), t \in \mathbb{R}^+)$. An analogue of Chung's law of the iterated logarithm is studied for $B^H$ and used to obtain the pointwise Hölder exponent of the local time. A kind of local asymptotic self-similarity is proved to be satisfied by the local time of $B^H$.

## 1. Introduction

Multifractional Brownian motion (mBm), $B^H = (B^{H(t)}(t), t \in \mathbb{R}^+)$, is a Gaussian process which extends fractional Brownian motion (fBm) by allowing the Hurst function to vary with time. This provides a tool to model systems whose regularity evolves over time, such as internet traffic or images. Recently, Boufoussi *et al.* [6] have investigated, under mild regularity conditions on $H(\cdot)$, the existence of jointly continuous local times of mBm. Their objective was to explore the link between the pointwise regularity of the local time and the irregularity of the underling process. This effect of inverse regularity was observed by Berman [5] for the case of uniform regularity.

The first aim of this paper is to establish estimates for the local and uniform moduli of continuity for the local time of mBm. Upper bounds for the moduli of continuity of local times have been obtained by Kôno [16] for Gaussian processes with stationary increments and, more recently, Csörgö *et al.* [8] have proved upper bounds for the moduli of continuity of the maxima of local times for stationary Gaussian processes and Gaussian processes with stationary increments. However, contrary to classical results, the moduli of continuity obtained in this paper depend on the point at which the regularity is studied. This is natural, since mBm has a regularity evolving over time. The approach of the present paper seems to be useful for extending the results of Kôno [16] and Csörgö *et al.* [8] to Gaussian processes without stationary increments.

Chung's form of the law of iterated logarithm for mBm is established in Section 4. This result is used to prove that the pointwise Hölder exponent of its local time with respect to time is equal to $1 - H(t)$, uniformly in $x$.

The second main objective of the paper concerns the question of whether the sample path properties of mBm can be transferred to its local time. If $H(\cdot)$ is not constant, the mBm is no longer self-similar. However, it is proved by Lévy-Véhel and Peltier ([17], Proposition 5) that if $H$ is $\beta$-Hölder continuous and $\sup_{t \in \mathbb{R}^+} H(t) < \beta$, then a property called *local asymptotic self-similarity* remains, defined as follows:

$$\lim_{\rho \longrightarrow 0^+} law\left\{\frac{B(t+\rho u) - B(t)}{\rho^{H(t)}}, u \in \mathbb{R}\right\} = law\{B^{H(t)}(u), u \in \mathbb{R}\}, \qquad (1)$$

where $B^{H(t)}$ is fBm with Hurst parameter $H(t)$. It is the purpose of Section 5 to prove that the local time of mBm has a kind of local asymptotic self-similarity. Through this result, we obtain some local limit theorems corresponding to the mBm.

We will use $C, C_1, \ldots$ to denote unspecified positive finite constants which may not necessarily be the same at each occurrence.

## 2. Preliminaries

In this section, we present some notation and collect facts about multifractional Brownian motion and its local times.

### 2.1. Multifractional Brownian motion

mBm was introduced independently by Lévy-Véhel and Peltier [17] and Benassi *et al.* [3]. The definition due to Lévy-Véhel and Peltier [17] is based on the moving average representation of fBm, where the constant Hurst parameter $H$ is replaced by a function $H(t)$ as follows:

$$B^{H(t)}(t) = \frac{1}{\Gamma(H(t)+1/2)} \bigg(\int_{-\infty}^{0} [(t-u)^{H(t)-1/2} - (-u)^{H(t)-1/2}]W(\mathrm{d}u) \\ + \int_{0}^{t}(t-u)^{H(t)-1/2}W(\mathrm{d}u)\bigg), \qquad t \geq 0, \quad (2)$$

where $H(\cdot): [0, \infty) \longrightarrow [\mu, \nu] \subset (0, 1)$ is a Hölder continuous function and $W$ is the standard Brownian motion defined on $(-\infty, +\infty)$. Benassi *et al.* [3] defined mBm by means of the harmonizable representation of fBm as follows:

$$\widehat{B}^{H(t)}(t) = \int_{\mathbb{R}} \frac{\mathrm{e}^{\mathrm{i}t\xi} - 1}{|\xi|^{H(t)+1/2}} \, \mathrm{d}\widehat{W}(\mathrm{d}\xi), \qquad (3)$$

where $\widehat{W}(\xi)$ is the Fourier transform of the series representation of white noise with respect to an orthonormal basis of $L^2(\mathbb{R})$ (we refer to Cohen [7] for the precise definition and for the fact that $\widehat{B}^H$ is real-valued).

Various properties of mBm have already been investigated in the literature, related, for instance, to its pointwise and uniform Hölder regularity, as well as the local Hausdorff dimension of its sample paths. More precisely, it is known from Lévy-Véhel and Peltier ([17], Proposition 8) that with probability one, for each $t_0 \geq 0$, the Hölder exponent at $t_0$ of mBm is $H(t_0)$. Recall that the pointwise Hölder exponent of a stochastic process $X$ at $t_0$ is defined by

$$\alpha_X(t_0, \omega) = \sup\left\{\alpha > 0, \lim_{\rho \to 0} \frac{X(t_0 + \rho, \omega) - X(t_0, \omega)}{\rho^\alpha} = 0\right\}. \quad (4)$$

In addition, according to the same authors, the local Hausdorff dimension of the graph is $2 - \min_{t \in [a,b]} H(t)$ almost surely for each interval $[a,b] \subset \mathbb{R}^+$.

According to the previous results, the regularity of mBm depends on the regularity of $H$. Furthermore, by using Lemma 3.1 in Boufoussi *et al.* [6], we can easily prove that the irregularity of $H$ implies the irregularity of mBm. More precisely, the points of discontinuity of $H$ are also discontinuities of $B$ (see Ayache [2], Proposition 1).

### 2.2. Local times

We end this section by briefly recalling some aspects of the theory of local times. For a comprehensive survey on local times of both random and non-random vector fields, we refer to Geman and Horowitz [12], Dozzi [10] and Xiao [21].

Let $X = (X(t), t \in \mathbb{R}^+)$ be a real-valued separable random process with Borel sample functions. For any Borel set $B \subset \mathbb{R}^+$, the occupation measure of $X$ on $B$ is defined as

$$\mu_B(A) = \lambda\{s \in B : X(s) \in A\}, \qquad \text{for all } A \in \mathcal{B}(\mathbb{R}),$$

where $\lambda$ is the one-dimensional Lebesgue measure on $\mathbb{R}^+$. If $\mu_B$ is absolutely continuous with respect to Lebesgue measure on $\mathbb{R}$, we say that $X$ has a local time on $B$ and define its local time, $L(B, \cdot)$, to be the Radon–Nikodym derivative of $\mu_B$. Here, $x$ is the so-called space variable and $B$ is the time variable.

By standard monotone class arguments, one can deduce that the local times have a measurable modification that satisfies the following *occupation density formula*: for every Borel set $B \subset \mathbb{R}^+$ and every measurable function $f : \mathbb{R} \to \mathbb{R}_+$,

$$\int_B f(X(t)) \, \mathrm{d}t = \int_\mathbb{R} f(x) L(B, x) \, \mathrm{d}x.$$

Recently, Boufoussi *et al.* [6] have proved that if the Hurst function is Hölder continuous with exponent $\beta$ and if $\sup_{t \geq 0} H(t) < \beta$, then the mBm has a jointly continuous local time, that is, the mapping $(t, x) \to L(t, x)$ is continuous. In addition, this local time has the following Hölder continuities. It satisfies, for any compact $U \subset \mathbb{R}$,

(i)

$$\sup_{x \in U} \frac{|L(t+h, x) - L(t, x)|}{|h|^\gamma} < +\infty \quad \text{a.s.,} \tag{5}$$

where $\gamma < 1 - H(t)$ and $|h| < \eta$, $\eta$ being a small random variable almost surely positive and finite;

(ii) for any $I \subset [0, T]$ with small length,

$$\sup_{x, y \in U, x \neq y} \frac{|L(I, x) - L(I, y)|}{|x - y|^\alpha} < +\infty \quad \text{a.s.,} \tag{6}$$

where $\alpha < (\frac{1}{2\sup_I H(t)} - \frac{1}{2}) \wedge 1$.

These results have been used to obtain the local and pointwise Hausdorff dimension of the level sets of mBm. We refer to Boufoussi *et al.* [6] for definitions and results.

## 3. Moduli of continuity of the local time

Throughout this section, the Hurst function $H: \mathbb{R}^+ \to [\mu, \nu] \subset (0, 1)$ is assumed to be measurable. The notation $B = (B(t), t \geq 0)$ means that either representation, (2) or (3), can be chosen. Moreover, we say that $H$ satisfies the condition $(\mathcal{H}_\beta)$ if

$$H \text{ is } \beta\text{-Hölder continuous with } \sup_{s \geq 0} H(s) < \beta.$$

We give the following improvement of Theorem 3.1 of Boufoussi *et al.* [6]. Note that the first part reproduces this theorem without assuming $(\mathcal{H}_\beta)$.

**Theorem 3.1.** *Consider a measurable function $H(\cdot): \mathbb{R}_+ \to [\mu, \nu] \subset (0, 1)$. The mBm with Hurst function $H(\cdot)$ given by (2) admits on any interval $[a, b] \subset [0, \infty)$ a square-integrable local time. Moreover, for both representations (2) and (3), if $(\mathcal{H}_\beta)$ holds, the existence of square-integrable local times implies that $H(t) < 1$ for almost all $t$.*

**Proof.** Let us write, for simplicity,

$$B^{H(t)}(t) = \int_{-\infty}^{t} K_{H(t)}(t, u) \, dW(u),$$

where

$$K_{H(t)}(t, u) = \frac{1}{\Gamma(H(t) + 1/2)} [(t - u)_+^{H(t) - 1/2} - (-u)_+^{H(t) - 1/2}].$$

For any $t > s$, taking the variance of $B^{H(t)}(t) - B^{H(s)}(s)$, we obtain

$$\text{Var}(B^{H(t)}(t) - B^{H(s)}(s)) \geq \text{Var}(B^{H(t)}(t) - B^{H(s)}(s)/W(u), u \leq s)$$
$$= \text{Var}(B^{H(t)}(t)/W(u), u \leq s), \tag{7}$$

where the last equality follows from the fact that $B^{H(s)}(s)$ is measurable with respect to $\sigma(W(u), u \leq s)$. Moreover, we can write

$$B^{H(t)}(t) = \int_{-\infty}^{t} K_{H(t)}(t,u)\,dW(u) = \int_{-\infty}^{s} K_{H(t)}(t,u)\,dW(u) + \int_{s}^{t} K_{H(t)}(t,u)\,dW(u).$$

Hence, by using the measurability of $\int_{-\infty}^{s} K_{H(t)}(t,u)\,dW(u)$ with respect to $\sigma(W(u)/u \leq s)$, we obtain

$$\mathrm{Var}(B^{H(t)}(t)/W(u), u \leq s) = \mathrm{Var}\left(\int_{s}^{t} K_{H(t)}(t,u)\,dW(u)/W(u), u \leq s\right)$$

$$= \mathrm{Var}\left(\int_{s}^{t} K_{H(t)}(t,u)\,dW(u)\right), \quad (8)$$

where, to obtain the last equality, we have used the fact that $\int_{s}^{t} K_{H(t)}(t,u)\,dW(u)$ is independent of $\sigma(W(u), u \leq s)$ (by the independence of the increments of the Brownian motion). Combining (7) and (8) and denoting $C = \sup_{u \in [\mu,\nu]} \Gamma(1/2 + u)$, we obtain

$$\mathrm{Var}(B^{H(t)}(t) - B^{H(s)}(s)) \geq \frac{1}{\Gamma(1/2 + H(t))^2} \int_{s}^{t} (t-r)^{2H(t)-1}\,dr$$

$$\geq \frac{1}{2\nu C^2}(t-s)^{2H(t)}.$$

Therefore,

$$\int_{[a,b]}\int_{[a,b]} (\mathrm{E}[B^{H(t)}(t) - B^{H(s)}(s)]^2)^{-1/2}\,ds\,dt$$

$$\leq \sqrt{2\nu}C \int_{[a,b]}\int_{[a,b]} |t-s|^{-\sup_{r \geq 0} H(r)}\,ds\,dt. \quad (9)$$

The last integral is finite because $\sup_{r \geq 0} H(r) < 1$. By Theorem 22.1 in Geman and Horowitz [12], $B^H$ has local time $L([a,b], \cdot) \in L^2(\mathbb{R})$.

We now prove the second point. According to Boufoussi *et al.* [6], we have, for any interval $[a,b] \subset [0, \infty)$ with small length,

$$E(B(t) - B(s))^2 \leq C_{\mu,\nu}|t-s|^{2H(t)}, \qquad \text{for } s, t \in [a,b]. \quad (10)$$

On the other hand, since the local time exists on the interval $[a,b]$, then according to Geman and Horowitz ([12], Theorem 22.1 expression (22.3)), the following integral is finite

$$\int_{a}^{b}\int_{a}^{b} \frac{1}{[E(B(t) - B(s))^2]^{1/2}}\,ds\,dt.$$

Then, (10) implies that

$$\int_a^b \int_a^b \frac{1}{(t-s)^{H(t)}} \, ds \, dt < \infty.$$

Consequently, $H(t) < 1$ for almost all $t \in [a, b]$. Since $\mathbb{R}_+$ is a countable union of small intervals, the result is proved. □

In what follows, we are interested in the local and pointwise oscillations (at each $t$) of the local time of mBm.

**Theorem 3.2.** *Let $\{B^{H(t)}(t), t \geq 0\}$ be the mBm given in (2) with Hurst function $H(\cdot)$ satisfying the assumption $(\mathcal{H}_\beta)$. Then, for every $t \in \mathbb{R}^+$ and any $x \in \mathbb{R}$, there exist positive and finite constants $C_1$ and $C_2$ such that*

$$\limsup_{\delta \to 0} \frac{L(t+\delta, B^{H(t)}(t)) - L(t, B^{H(t)}(t))}{\delta^{1-H(t)} (\log \log(\delta^{-1}))^{H(t)}} \leq C_1 \qquad a.s., \tag{11}$$

$$\limsup_{\delta \to 0} \frac{L(t+\delta, x) - L(t, x)}{\delta^{1-H(t)} (\log \log(\delta^{-1}))^{H(t)}} \leq C_2 \qquad a.s. \tag{12}$$

**Proof.** Let $t \geq 0$ be a fixed point, suppose $0 < h < 1$ and define $H_{t,t+h} = \sup_{s \in [t,t+h]} H(s)$. According to the Fourier analytic approach of Berman [4] (see also Davies [9], expression (27)), we have

$$E[L(t+h, B^{H(t)}(t)) - L(t, B^{H(t)}(t))]^m$$
$$= \frac{1}{(2\pi)^m} \int_{[t,t+h]^m} \int_{\mathbb{R}^m} E(e^{i \sum_{j=1}^m u_j (B^{H(s_j)}(s_j) - B^{H(t)}(t))}) \prod_{j=1}^m du_j \prod_{j=1}^m ds_j.$$

Let $\widetilde{B}(s) = B^{H(s)}(s) - B^{H(t)}(t)$, $s \geq 0$, and denote by $R(s_1, \ldots, s_m)$ the covariance matrix of $(\widetilde{B}(s_1), \ldots, \widetilde{B}(s_m))$ for distinct $s_1, \ldots, s_m$. Let $U = (u_1, \ldots, u_m) \in \mathbb{R}^m$ and let $U^T$ denote the transpose of $U$. According to (16) below, we have $\det R(s_1, \ldots, s_m) > 0$. Hence, the change of variable $V = R^{1/2} U$ implies that

$$\int_{\mathbb{R}^m} E(e^{i \sum_{j=1}^m u_j \widetilde{B}(s_j)}) \, du_1 \cdots du_m = \frac{(2\pi)^{m/2}}{(\det R(s_1, \ldots, s_m))^{1/2}}.$$

Therefore,

$$E[L(t+h, B^{H(t)}(t)) - L(t, B^{H(t)}(t))]^m$$
$$= \frac{1}{(2\pi)^{m/2}} \int_{[t,t+h]^m} \frac{1}{(\det R(s_1, \ldots, s_m))^{1/2}} \, ds_1 \cdots ds_m$$
$$= \frac{m!}{(2\pi)^{m/2}} \int_{t < s_1 < \cdots < s_m < t+h} \frac{1}{(\det R(s_1, \ldots, s_m))^{1/2}} \, ds_1 \cdots ds_m. \tag{13}$$

In addition,

$$\det R(s_1,\ldots,s_m)$$
$$= \operatorname{Var}(\widetilde{B}(s_1))\operatorname{Var}(\widetilde{B}(s_2)/\widetilde{B}(s_1))\cdots\operatorname{Var}(\widetilde{B}(s_m)/\widetilde{B}(s_1),\ldots,\widetilde{B}(s_{m-1})). \quad (14)$$

Using arguments similar to those used in the proof of Theorem 3.1, we obtain, for any $r, s \in [t, t+h]$ such that $r < s$,

$$\operatorname{Var}(\widetilde{B}(s)/\widetilde{B}(u), u \in A, u \leq r) \geq \operatorname{Var}(\widetilde{B}(s) - \widetilde{B}(r)/W(u), u \leq r)$$
$$\geq \operatorname{Var}(B^{H(s)}(s)/W(u), u \leq r)$$
$$\geq \frac{1}{2\nu C^2}(s-r)^{2H(s)}, \quad (15)$$

where $C = \sup_{x \in [\mu,\nu]} \Gamma(x+1/2)$. Combining (14) and (15), we obtain

$$\det R(s_1,\ldots,s_m) \geq \frac{1}{(2\nu C^2)^m} \prod_{j=1}^{m} (s_j - s_{j-1})^{2H(s_j)}$$
$$\geq \frac{1}{(2\nu C^2)^m} \prod_{j=1}^{m} (s_j - s_{j-1})^{2H_{t,t+h}}, \quad (16)$$

where $s_0 = 0$ and $(s_j - s_{j-1})^{H(s_j)} \geq (s_j - s_{j-1})^{H_{t,t+h}}$, since $(s_j - s_{j-1}) < 1$.

According to (13), we have

$$E[L(t+h, B^{H(t)}(t)) - L(t, B^{H(t)}(t))]^m$$
$$\leq \left(\frac{\nu C^2}{\pi}\right)^{m/2} m! \int_{t<s_1<\cdots<s_m<t+h} \prod_{j=1}^{m} \frac{1}{(s_j - s_{j-1})^{H_{t,t+h}}} \, ds_1 \cdots ds_m.$$

Now, by an elementary calculation (cf. Ehm [11]), for all $m \geq 1$, $h > 0$ and $b_j < 1$,

$$\int_{t<s_1<\cdots<s_m<t+h} \prod_{j=1}^{m} (s_j - s_{j-1})^{-b_j} \, ds_1 \cdots ds_m = h^{m - \sum_{j=1}^{m} b_j} \frac{\prod_{j=1}^{m} \Gamma(1-b_j)}{\Gamma(1+m-\sum_{j=1}^{m} b_j)}.$$

Therefore,

$$E[L(t+h, B^{H(t)}(t)) - L(t, B^{H(t)}(t))]^m$$
$$\leq m! \left(\frac{\nu C^2}{\pi}\right)^{m/2} h^{m(1-H_{t,t+h})} \frac{(\Gamma(1-H_{t,t+h}))^m}{\Gamma(1+m(1-H_{t,t+h}))}.$$

According to Stirling's formula, we have $(m!/\Gamma(1+m(1-H_{t,t+h}))) \leq M^m m!^{H_{t,t+h}}, m \geq 2$, for a suitable finite number $M$. Therefore,

$$E\left[\frac{L(t+h, B^{H(t)}(t)) - L(t, B^{H(t)}(t))}{h^{1-H_{t,t+h}}}\right]^m \leq C^m m!^{H_{t,t+h}}. \tag{17}$$

We shall now prove that for any $K > 0$, there exists a positive and finite constant $A > 0$, depending on $t$, such that for sufficiently small $u$,

$$P\left(L(t+h, B^{H(t)}(t)) - L(t, B^{H(t)}(t)) \geq \frac{Ah^{1-H_{t,t+h}}}{u^{H_{t,t+h}}}\right) \leq \exp(-K/u). \tag{18}$$

First consider $u$ of the form $u = 1/m$. Combining Chebyshev's inequality and (17), we obtain

$$P\left(\frac{L(t+h, B^{H(t)}(t)) - L(t, B^{H(t)}(t))}{h^{1-H_{t,t+h}}} \geq Am^{H_{t,t+h}}\right)$$

$$\leq E\left[\frac{L(t+h, B^{H(t)}(t)) - L(t, B^{H(t)}(t))}{Ah^{(1-H_{t,t+h})}m^{H_{t,t+h}}}\right]^m$$

$$\leq \frac{C^m}{A^m}\left(\frac{1}{m}\right)^{mH_{t,t+h}} (m!)^{H_{t,t+h}}.$$

Again using Stirling's formula, the last expression is at most $\frac{C^m}{A^m}(2\pi m)^{H_{t,t+h}/2} e^{-H_{t,t+h}m}$. This can be written as

$$\exp\left(m[\log(C/A) - H_{t,t+h}] + \frac{H_{t,t+h}}{2}[\log(m) + \log(2\pi)]\right). \tag{19}$$

Choose $A > C$ and $m_0$ large such that for any $m \geq m_0$, to dominate (19) by $e^{-2Km}$. Moreover, for $u$ sufficiently small, there exists $m \geq m_0$ such that $u_{m+1} < u < u_m$ and since $m \geq 1$, $\frac{m}{m+1} \geq \frac{1}{2}$. This proves (18).

In addition, if we take $u(h) = 1/\log\log(1/h)$ and first consider $h_m$ of the form $2^{-m}$, then (18) implies

$$P(L(t+h_m, B^{H(t)}(t)) - L(t, B^{H(t)}(t)) \geq Ah_m^{1-H_{t,t+h_m}}(\log\log(1/h_m))^{H_{t,t+h_m}}) \leq m^{-2}$$

for large $m$. Consequently, by using the Borel–Cantelli lemma and monotonicity arguments, we obtain

$$\frac{L(t+h, B^{H(t)}(t)) - L(t, B^{H(t)}(t))}{h^{1-H_{t,t+h}}} \leq A(\log\log(1/h))^{H_{t,t+h}} \quad \text{a.s.}$$

In addition, the Hölder regularity of $H$ implies that

$$\frac{\log\log(1/h)^{H_{t,t+h}}}{\log\log(1/h)^{H(t)}} = e^{Kh^\beta \log\log\log 1/h} \to 1 \quad \text{as } h \to 0^+.$$

Therefore,
$$(\log\log(1/h))^{H_{t,t+h}} \leq K(\log\log(1/h))^{H(t)}$$

for sufficiently small $h$. Moreover, using the same arguments (see Remark 3.7 in Boufoussi *et al.* [6]), we obtain

$$h^{1-H_{t,t+h}} \leq Kh^{1-H(t)}. \tag{20}$$

This completes the proof of (11). Since (12) is proved in the same manner, the proof is omitted here. $\square$

**Remark 3.3.** The result of the previous theorem, with $\sup_{s\in[t,t+h]} H(s)$ instead of $H(t)$, may be proven even when $(\mathcal{H}_\beta)$ is not satisfied.

We shall also prove the following uniform result.

**Theorem 3.4.** *Let $\{B^{H(t)}(t), t \in [0,1]\}$ be the mBm in (2) with arbitrary Hurst function and define $H_{0,1} = \sup_{s\in[0,1]} H(s)$. Then, for all $x \in \mathbb{R}$, there exists $C_3 > 0$ such that*

$$\limsup_{\substack{|t-s|\searrow 0^+ \\ t,s\in[0,1]}} \frac{|L(t,x)-L(s,x)|}{|t-s|^{1-H_{0,1}}(\log(1/|t-s|))^{H_{0,1}}} \leq C_3 \qquad a.s. \tag{21}$$

**Proof.** By using the same arguments as above, without using the Hölder continuity of $H$, we show that with probability one, there exists $m_0(\omega)$ such that for all $m \geq m_0(\omega)$ and $k = 1, \ldots, 2^m$, we have

$$L(k2^{-m},x) - L((k-1)2^{-m},x) \leq A2^{-(1-H_{0,1})m}(m\log(2))^{H_{0,1}} \qquad \text{a.s.}$$

By then proceeding essentially as in Kôno [16], we prove (21). $\square$

**Remark 3.5.** The results in this section are based on the following property of the mBm:

$$\operatorname{Var}(B^{H(t+h)}(t+h)/B^{H(s)}(s) : 0 \leq s \leq t) \geq Ch^{2H_{t,t+h}}, \tag{22}$$

for sufficiently small $h$. If $H$ is constant, this property is called *one-sided strong local nondeterminism*; we refer to Monrad and Rootzén [20], expression (2.2), and the references therein for definition and applications. Property (22) is satisfied by many other processes with multifractal behavior, such as the Riemann–Liouville mBm $\{X(t), t \geq 0\}$ introduced in Lim [18] by replacing, in the definition of the Riemann–Liouville fBm, the constant Hurst parameter $H$ by a regular function as follows:

$$X(t) = \frac{1}{\Gamma(H(t)+1/2)} \int_0^t (t-u)^{H(t)-1/2} \, dW_u.$$

Thus our results extend those of Kôno [16] to Gaussian processes without the stationary increments property.

*Remark 3.6.* A closer examination of the proof of Theorem 2.1 in Csörgö *et al.* [8] shows that the property of stationarity of increments is only used to prove two main ingredients presented in their Lemmas 3.1 and 3.3 and which can, in our case, be easily replaced by the following properties:

1. Let $\mathbb{A}_n$ be the covariance matrix of a Gaussian vector $\{\zeta_i; 1 \leq i \leq n\}$. It is then known that the conditional variance can be written as

$$\mathrm{Var}(\zeta_i/\zeta_l, l \neq i, 1 \leq l \leq n) = \frac{|\mathbb{A}_n|}{|\mathbb{A}_n^{(i)}|},$$

where $\mathbb{A}_n^{(i)}$ is the submatrix of $\mathbb{A}_n$ obtained by deleting the $i$th row and column.

2. If a stochastic process has the property (22), with $H$ constant, then, according to (14), its covariance matrix satisfies

$$\det R(t_1, \ldots, t_m) \geq C^m \prod_{j=1}^{m}(t_j - t_{j-1})^{2H}, \qquad \text{for } t < t_1 < \cdots < t_m < t+h \text{ and } t_0 = 0.$$

The arguments of Csörgö *et al.* ([8], Theorem 2.1) can then be adapted to prove that the class of processes satisfying (22) satisfies

$$\limsup_{T \to \infty} \sup_{0 \leq t \leq b_T} \frac{L(t+a_T, x) - L(t, x)}{a_T^{1-H}(\log(b_T/a_T) + \log\log(a_T + 1/a_T))^H} < C < \infty \qquad \text{a.s.,} \qquad (23)$$

where $a_T$ and $b_T$ are non-negative functions in $T \geq 0$ such that $\frac{1+b_T}{a_T} \longrightarrow +\infty$ as $T \to +\infty$.

## 4. Chung's law of the iterated logarithm

The following theorem establishes an analogue of Chung's law of the iterated logarithm for mBm. The result for fBm has been obtained by Monrad and Rootzén [20].

**Theorem 4.1.** *Let $B$ be a multifractional Brownian motion and assume that $H$ satisfies condition $(\mathcal{H}_\beta)$. Then the following Chung-type law of iterated logarithm holds:*

$$\liminf_{\delta \to 0} \sup_{s \in [t_0, t_0 + \delta]} \frac{|B(s) - B(t_0)|}{(\delta/\log|\log(\delta)|)^{H(t_0)}} = C_{H(t_0)} \qquad a.s., \qquad (24)$$

*where the constant $C_{H(t_0)}$ is the one appearing in Chung's law for the fBm of Hurst parameter $H(t_0)$.*

**Proof.** We first introduce the process $\{\widetilde{B}(t) = B(t) - B^{H(t_0)}(t), t \geq 0\}$, where $B^{H(t_0)}$ denotes a fBm with Hurst parameter $H(t_0)$. According to Theorem 3.3 of Monrad and Rootzén [20], the fBm with Hurst exponent $H(t_0)$ satisfies (24). Furthermore,

$$B(t) - B(t_0) = (B(t) - B^{H(t_0)}(t)) + (B^{H(t_0)}(t) - B^{H(t_0)}(t_0)),$$
$$= \widetilde{B}(t) + (B^{H(t_0)}(t) - B^{H(t_0)}(t_0)).$$

(24) will then be proven if we show that

$$\lim_{\delta \to 0} \sup_{t \in [t_0, t_0+\delta]} \frac{|\widetilde{B}(t)|}{(\delta/\log|\log(\delta)|)^{H(t_0)}} = 0 \quad \text{a.s.}$$

On the other hand, Lemma 3.1 in Boufoussi *et al.* [6] implies that there exists a positive constant $K$ such that

$$\sup_{t \in [t_0, t_0+\delta]} E(\widetilde{B}(t))^2 \leq K\delta^{2\beta}.$$

Hence, according to Theorem 2.1 in Adler ([1], page 43) and a symmetry argument, we have

$$P\left(\sup_{t \in [t_0, t_0+\delta]} |\widetilde{B}(t)| \geq u\right) \leq 2P\left(\sup_{t \in [t_0, t_0+\delta]} \widetilde{B}(t) \geq u\right)$$
$$\leq 4\exp\left(-\frac{(u - E(\sup_{t \in [t_0, t_0+\delta]} \widetilde{B}(t)))^2}{K\delta^{2\beta}}\right). \quad (25)$$

For the sake of simplicity, let $\Lambda = \sup_{t \in [t_0, t_0+\delta]} \widetilde{B}(t)$. By (25), we obtain

$$E(\Lambda) \leq \int_0^{+\infty} P\left(\sup_{t \in [t_0, t_0+\delta]} |\widetilde{B}(t)| > x\right) dx$$
$$\leq 4 \int_0^{+\infty} \exp\left(-\frac{[x - E\Lambda]^2}{K\delta^{2\beta}}\right) dx$$
$$= \frac{4\sqrt{K}\delta^\beta}{\sqrt{2}} \int_{-\sqrt{2/K}E\Lambda/\delta^\beta}^{\infty} e^{-y^2/2} dy$$
$$\leq 4\sqrt{K\pi}\delta^\beta.$$

It follows that

$$(u - E\Lambda)^2 \geq \tfrac{1}{2}u^2 - (E\Lambda)^2 \geq \tfrac{1}{2}u^2 - 16K\pi\delta^{2\beta}.$$

Consequently, (25) becomes

$$P\left(\sup_{t \in [t_0, t_0+\delta]} |\widetilde{B}(t)| \geq u\right) \leq C\exp\left(-\frac{u^2}{K\delta^{2\beta}}\right). \quad (26)$$

Since $H(t_0) < \beta$, there exists $0 < \xi < \beta - H(t_0)$. Consider $\delta_n = n^{1/(2(\xi+H(t_0)-\beta))}$ and $u_n = \delta_n^{H(t_0)+\xi}$. Therefore, according to (26), we have

$$\sum_{n=1}^{\infty} P\left(\sup_{t \in [t_0, t_0+\delta_n]} |\widetilde{B}(t)| \geq u_n\right) \leq \sum_{n=1}^{\infty} \exp\left(-\frac{1}{K}n\right) < \infty.$$

It follows from the Borel–Cantelli lemma that there exists $n_0 = n_0(\omega)$ such that for all $n \geq n_0$, $\sup_{s \in [t_0, t_0+\delta_n]} |\widetilde{B}(s)| \leq \delta_n^{H(t_0)+\xi}$ almost surely. Furthermore, for $\delta_{n+1} \leq \delta \leq \delta_n$, we have

$$\sup_{s \in [t_0, t_0+\delta]} |\widetilde{B}(s)| \leq \sup_{s \in [t_0, t_0+\delta_n]} |\widetilde{B}(s)|$$

$$\leq \delta_n^{H(t_0)+\xi} \leq \delta^{H(t_0)+\xi}\left(\frac{\delta_n}{\delta_{n+1}}\right)^{H(t_0)+\xi}$$

$$\leq 2^{\theta} \delta^{H(t_0)+\xi} \quad \text{a.s.,}$$

where $\theta = \frac{H(t_0)+\xi}{2(\beta-H(t_0)-\xi)}$. Hence,

$$\lim_{\delta \to 0} \sup_{t \in [t_0, t_0+\delta]} \frac{|\widetilde{B}(t)|}{(\delta/\log|\log(\delta)|)^{H(t_0)}} \leq \lim_{\delta \to 0} \delta^{\xi}(\log|\log(\delta)|)^{H(t_0)} = 0 \quad \text{a.s.}$$

This completes the proof of the theorem. □

*Remark 4.2.* Observe that by the ideas used to prove the previous theorem, many laws of iterated logarithm (LIL) proved for fBm can now be obtained for mBm and with the same constants. For example, we have the LIL (cf. Li and Shao [19], equation (7.5) for the fBm)

$$\limsup_{\delta \to 0} \sup_{s \in [t_0, t_0+\delta]} \frac{|B(s) - B(t_0)|}{\delta^{H(t_0)}(\log|\log(\delta)|)^{1/2}} = \sqrt{2} V_{H(t_0)}(B) \quad \text{a.s.,} \quad (27)$$

where

$$V_{H(t_0)}(\widehat{B}) = \sqrt{\frac{\pi}{H(t_0)\Gamma(2H(t_0))\sin(\pi H(t_0))}}$$

and

$$V_{H(t_0)}(B^H) = \frac{(\int_{-\infty}^{0}[(1-s)^{H(t_0)-1/2} - (-s)^{H(t_0)-1/2}]^2 \, ds)^{1/2} + 1/(2H(t_0))}{\Gamma(H(t_0)+1/2)}.$$

Other laws of iterated logarithm and uniform moduli of continuity have been obtained for $\widehat{B}$ by Benassi *et al.* [3], Theorem 1.7, via wavelet techniques.

Since the local time is a two-parameter process, we make precise the meaning of the pointwise Hlder exponent in time, uniformly in space. The pointwise Hlder exponent of $\sup_{x\in\mathbb{R}} L(\cdot,x)$ is defined by

$$\alpha_L(t) = \sup\left\{\alpha > 0, \limsup_{\delta\to 0} \sup_{x\in\mathbb{R}} \frac{L(t_0+\delta,x) - L(t_0,x)}{\delta^\alpha} = 0\right\}. \tag{28}$$

Since the fraction above is positive, the lim sup may be replaced by lim or lim inf.

We are now able to prove the following.

**Proposition 4.3.** *Let $B$ be a multifractional Brownian motion and assume that $H$ satisfies the condition $(\mathcal{H}_\beta)$. The pointwise Hölder exponent $\alpha_L$ of $\sup_{x\in\mathbb{R}} L(\cdot,x)$ at $t$ satisfies*

$$\alpha_L(t) = 1 - H(t) \qquad a.s.$$

**Proof.** The lower bound was already given in Boufoussi *et al.* [6], Corollary 3.6. The upper bound is a consequence of Chung's law of the mBm. Indeed, since the local time vanishes outside the range of $B$, we obtain

$$\begin{aligned}
\delta &= \int_{\mathbb{R}} L([t_0, t_0+\delta], x)\,\mathrm{d}x \\
&\leq \sup_{x\in\mathbb{R}} L([t_0, t_0+\delta], x) \sup_{s,t\in[t_0,t_0+\delta]} |B(s) - B(t)| \\
&\leq 2 \sup_{x\in\mathbb{R}} L([t_0, t_0+\delta], x) \sup_{s\in[t_0,t_0+\delta]} |B(s) - B(t_0)|.
\end{aligned} \tag{29}$$

Combining (24) and (29), we obtain that there exists a positive constant $C$ such that

$$\limsup_{\delta\to 0} \sup_{x\in\mathbb{R}} \frac{L(t_0+\delta,x) - L(t_0,x)}{\delta^{1-H(t_0)}(\log\log(\delta^{-1}))^{H(t_0)}} \geq C \qquad \text{a.s.},$$

which, together with the definition of the pointwise Hölder exponent, proves the result. □

## 5. Asymptotic results

It is well known that techniques for proving limit theorems related to self-similar processes use the self-similarity of their local times. It is natural to expect the same when dealing with locally asymptotically self-similar processes (LASS for short). Thus it will be of some interest to know if the local times satisfy a kind of LASS property.

### 5.1. LASS for local times

The answer to the preceding question is affirmative in the case of the mBm and the result is given by the following theorem.

**Theorem 5.1.** *Let $B$ be a multifractional Brownian motion and assume that $H$ satisfies the condition $(\mathcal{H}_\beta)$. Then, for any fixed $t_0$, the local time of mBm is locally asymptotically self-similar with parameter $1 - H(t_0)$, in the sense that for every $x \in \mathbb{R}$, the processes $\{Y_\rho(t,x), t \in [0,1]\}_{\rho > 0}$, defined by*

$$Y_\rho(t,x) = \frac{L(t_0 + \rho t, \rho^{H(t_0)}x + B(t_0)) - L(t_0, \rho^{H(t_0)}x + B(t_0))}{\rho^{1-H(t_0)}},$$

*converge in law to the local time, $\{\ell(t,x), t \in [0,1]\}$, of the fBm $B^{H(t_0)}$ with Hurst parameter $H(t_0)$, that is,*

$$\lim_{\rho \to 0} law\{Y_\rho(t,x), t \in [0,1]\} = law\{\ell(t,x), t \in [0,1]\}, \tag{30}$$

*where the convergence is in the space of continuous functions endowed with the norm of uniform convergence.*

The motivation for considering the processes $Y_\rho$ follows from the occupation density formula and the LASS property (1) of mBm. Indeed, according to the occupation density formula, the local time has the following representation:

$$L(t,x) = \frac{1}{2\pi} \int_{-\infty}^{+\infty} e^{-iux} \left( \int_0^t e^{iuB(s)} \, ds \right) du.$$

Consequently,

$$Y_\rho(t,x) = \frac{1}{2\pi \rho^{1-H(t_0)}} \int_\mathbb{R} e^{-iy(\rho^{H(t_0)}x + B(t_0))} \int_{t_0}^{t_0 + \rho t} e^{iyB(s)} \, ds \, dy.$$

Using the change of variables $r = \frac{s-t_0}{\rho}$ and $\rho^{H(t_0)}y = v$, the right-hand side of the previous expression becomes

$$\frac{1}{2\pi} \int_\mathbb{R} e^{-ivx} \int_0^t \exp\left(iv \frac{B(t_0 + \rho r) - B(t_0)}{\rho^{H(t_0)}}\right) dr \, dv,$$

which is the local time of the Gaussian process $\{B^\rho(r) = \frac{B(t_0 + \rho r) - B(t_0)}{\rho^{H(t_0)}}, r \in [0,1]\}$. We need the following lemma for the proof of finite-dimensional convergence.

**Lemma 5.1.** *Let $\{X(t), t \in [0,1]\}$ be a stochastic process in the Skorohod space $D([0,1])$ and define, for a fixed interval $I = [a,b] \subset \mathbb{R}$, the map*

$$\phi_I(X)(t) = \int_0^t 1_{\{X(s) \in I\}} \, ds.$$

*If $\{X_n(t), t \in [0,1]\}_{n \geq 1}$ is a family of processes which converges in law in $D([0,1])$ to $X$, then $\phi_I(X_n)$ converges in law to $\phi_I(X)$.*

**Proof.** The lemma is a consequence of the continuity of the map $X(\cdot) \to \phi_I(X(\cdot))$ in $J_1$ topology on $D([0,1])$ at almost all sample points of the process X, which is proved on page 11 in Kesten and Spitzer [15]. □

**Proof of Theorem 5.1.** To prove the convergence in law, we proceed in two steps. First, we prove the tightness of the family $\{Y_\rho(t,x), t \in [0,1]\}_{\rho>0}$ in the space of continuous functions. By using (17) and (20), for sufficiently small $\rho$, we obtain

$$E|Y_\rho(t,x) - Y_\rho(s,x)|^m$$
$$= \frac{E[L(t_0 + \rho t, \rho^{H(t_0)}x + B(t_0)) - L(t_0 + \rho s, \rho^{H(t_0)}x + B(t_0))]^m}{\rho^{(1-H(t_0))m}}$$
$$\leq C_m |t-s|^{(1-H(t_0))m}.$$

We can take $m > \frac{1}{1-H(t_0)}$ to prove the tightness.

Now, we prove the convergence of the finite-dimensional distributions of $Y_\rho$, as $\rho$ tends to 0, to those of the local time $\ell$ of the fBm $B^{H(t_0)}$ with Hurst parameter $H(t_0)$. We need to show that for any $d \geq 1$, $a_1, \ldots, a_d \in \mathbb{R}$ and $t_1, \ldots, t_d \in [0,1]$, the following convergence holds:

$$\sum_{j=1}^d a_j Y_\rho(t_j, x) \xrightarrow{\mathcal{W}} \sum_{j=1}^d a_j \ell(t_j, x) \qquad \text{as } \rho \to 0.$$

We will show the convergence of the corresponding characteristic function. More precisely, we will prove that

$$\left| E \exp\left[i\lambda \sum_{j=1}^d a_j Y_\rho(t_j, x)\right] - E \exp\left[i\lambda \sum_{j=1}^d a_j \ell(t_j, x)\right] \right| \longrightarrow 0 \qquad \text{as } \rho \to 0.$$

We introduce the following notation:

$$\phi_{\varepsilon,x}(X)(t) = \frac{1}{\varepsilon} \int_0^t 1_{[x, x+\varepsilon]}(X(s)) \, ds,$$

$$I_1^{\varepsilon,\rho} = \left| E \exp\left[i\lambda \sum_{j=1}^d a_j Y_\rho(t_j, x)\right] - E \exp\left[i\lambda \sum_{j=1}^d a_j \phi_{\varepsilon,x}(B^\rho)(t_j)\right] \right|,$$

$$I_2^{\varepsilon,\rho} = \left| E \exp\left[i\lambda \sum_{j=1}^d a_j \phi_{\varepsilon,x}(B^{H(t_0)})(t_j)\right] - E \exp\left[i\lambda \sum_{j=1}^d a_j \phi_{\varepsilon,x}(B^\rho)(t_j)\right] \right|$$

and

$$I_3^\varepsilon = \left| E \exp\left[i\lambda \sum_{j=1}^d a_j \phi_{\varepsilon,x}(B^{H(t_0)})(t_j)\right] - E \exp\left[i\lambda \sum_{j=1}^d a_j \ell(t_j, x)\right] \right|.$$

Therefore,

$$\left| \mathrm{E} \exp\left[i\lambda \sum_{j=1}^{d} a_j Y_\rho(t_j, x)\right] - \mathrm{E} \exp\left[i\lambda \sum_{j=1}^{d} a_j \ell(t_j, x)\right] \right| \leq I_1^{\varepsilon,\rho} + I_2^{\varepsilon,\rho} + I_3^{\varepsilon,\rho}. \quad (31)$$

On the other hand, $Y_\rho$ is the local time of $B^\rho$ and by using the mean value theorem and the occupation density formula, we obtain

$$I_1^{\varepsilon,\rho} \leq C \max_{1 \leq j \leq d} E|Y_\rho(t_j, x) - \phi_{\varepsilon,x}(B^\rho)(t_j)|$$

$$= C \max_{1 \leq j \leq d} E\left| \frac{1}{\varepsilon} \int_x^{x+\varepsilon} Y_\rho(t_j, y) \, \mathrm{d}y - Y_\rho(t_j, x) \right|. \quad (32)$$

Moreover, since the stochastic process $\{Y_\rho(t, y), y \in \mathbb{R}\}$ is almost surely continuous in $y$ for every $t$, according to the dominated convergence theorem, (32) converges to zero as $\varepsilon$ tends to zero independently of $\rho$.

We now deal with $I_2^{\varepsilon,\rho}$. Since the family of processes $\{B^\rho(t), t \in [0,1]\}_{\rho>0}$ converges in distribution to the fBm $\{B^{H(t_0)}(t), t \in [0,1]\}$ with Hurst parameter $H(t_0)$, the second term converges to zero as $\rho$ tends to 0 by Lemma 5.1.

The last term in (31) is treated in a similar way as the first and the proof of the finite-dimensional convergence is complete. $\square$

### 5.2. Limit theorems

The following result is an immediate consequence of the LASS property of mBm.

**Proposition 5.2.** *Let $B$ be a multifractional Brownian motion, assume that $H$ satisfies the condition $(\mathcal{H}_\beta)$ and denote by $\ell(t, x)$ the local time of the fBm with Hurst parameter $H(t_0)$. Then, for every $f \in L^1(\mathbb{R})$ which is locally Riemann integrable, with compact support and such that $\int_\mathbb{R} f(x) \, \mathrm{d}x \neq 0$, the following convergence in law holds:*

$$\lim_{\lambda \to \infty} \lim_{\rho \to 0^+} \frac{1}{\lambda^{1-H(t_0)}} \int_0^{\lambda t} f\left( \frac{B(\rho s + t_0) - B(t_0)}{\rho^{H(t_0)}} \right) \mathrm{d}s = \int_\mathbb{R} f(x) \, \mathrm{d}x \cdot \ell(t, 0).$$

If $H$ is constant, then $B$ is fBm. Being self-similar with stationary increments, the term in the left-hand side of the previous formula is identical in law to $\frac{1}{\lambda^{1-H}} \int_0^{\lambda t} f(B(s)) \, \mathrm{d}s$. We retrieve the result of Kasahara and Kosugi [14] for the fBm.

**Proof of Proposition 5.2.** The proof relies on the LASS property and the self-similarity of fBm. We sketch it for the sake of completeness.

Combining the fact that $f$ is locally Riemann integrable, the LASS property of mBm and Theorem VI.4.2 in Gihman and Skorohod [13], we obtain

$$\int_0^{\lambda t} f\left( \frac{B(\rho s + t_0) - B(t_0)}{\rho^{H(t_0)}} \right) \mathrm{d}s \xrightarrow{(d)} \int_0^{\lambda t} f(B^{H(t_0)}(s)) \, \mathrm{d}s \qquad \text{as } \rho \to 0^+.$$

Moreover, using the occupation density formula and the self-similarity of fBm, we obtain

$$\frac{1}{\lambda^{1-H(t_0)}} \int_0^{\lambda t} f(B^{H(t_0)}(s)) \, ds \stackrel{(d)}{=} \int_{\mathbb{R}} f(x) \ell(t, \lambda^{-H(t_0)} x) \, dx,$$

where $\stackrel{(d)}{=}$ denotes equality in distribution. The result of the theorem then follows from the continuity of $\ell(t, x, \omega)$ with respect to the space variable. $\square$

However, the presence of a double limit may be not convenient and the result may be better if $\lambda$ and $\rho$ are dependent. We prove the following local limit theorem.

**Theorem 5.3.** *Let B be a multifractional Brownian motion and assume that H satisfies the condition $(\mathcal{H}_\beta)$. Denote by $\ell(t,x)$ the local time of the fBm with Hurst parameter $H(t_0)$ and consider $f \in L^1(\mathbb{R})$ such that*

$$\int_{\mathbb{R}} |f(x)||x|^\xi \, dx < \infty \qquad \text{for some } 0 < \xi < \frac{1}{2\sup_{t \geq 0} H(t)} - \frac{1}{2} \quad \text{and} \quad \int_{\mathbb{R}} f(x) \, dx \neq 0. \quad (33)$$

*The following convergence in law then holds*

$$\frac{1}{\psi(\rho)} \int_{t_0}^{\rho t + t_0} f\left(\frac{B(s) - B(t_0) - \rho^{H(t_0)} y}{\theta(\rho)}\right) ds \xrightarrow{(d)} \int_{\mathbb{R}} f(x) \, dx \cdot \ell(t, y), \qquad \text{as } \rho \to 0^+,$$

*where $\theta(\cdot)$ and $\psi(\cdot)$ satisfy $\frac{\psi(\rho)}{\theta(\rho)} = \rho^{1-H(t_0)}$ and $\frac{\theta(\rho)}{\rho^{H(t_0)}} = o(1)$.*

**Proof.** Using the occupation density formula, we obtain

$$\frac{1}{\psi(\rho)} \int_{t_0}^{\rho t + t_0} f\left(\frac{B(s) - B(t_0) - \rho^{H(t_0)} y}{\theta(\rho)}\right) ds$$

$$= \int_{\mathbb{R}} f(x) \frac{L(\rho t + t_0, \theta(\rho) x + \rho^{H(t_0)} y + B(t_0)) - L(t_0, \theta(\rho) x + \rho^{H(t_0)} y + B(t_0))}{\rho^{1-H(t_0)}} \, dx$$

$$= \frac{L(\rho t + t_0, \rho^{H(t_0)} y + B(t_0)) - L(t_0, \rho^{H(t_0)} y + B(t_0))}{\rho^{1-H(t_0)}} \times \int_{\mathbb{R}} f(x) \, dx \qquad (34)$$

$$+ \int_{\mathbb{R}} f(x) \frac{L(I_0, \theta(\rho) x + \rho^{H(t_0)} y + B(t_0)) - L(I_0, \rho^{H(t_0)} y + B(t_0))}{\rho^{1-H(t_0)}} \, dx, \qquad (35)$$

where, in the last expression, we let $I_0 = [t_0, t_0 + \rho t]$ for simplicity. According to Theorem 5.1, the expression (34) converges in distribution to $\ell(t, y) \times \int_{\mathbb{R}} f(x) \, dx$. It now suffices to prove that (35) converges to 0 in some strong sense. We have

$$\mathbb{E}\left| \int_{\mathbb{R}} f(x) \frac{L(I_0, \theta(\rho) x + \rho^{H(t_0)} y + B(t_0)) - L(I_0, \rho^{H(t_0)} y + B(t_0))}{\rho^{1-H(t_0)}} \, dx \right|$$

$$\leq \int_{\mathbb{R}} |f(x)| \left\| \frac{L(I_0, \theta(\rho)x + \rho^{H(t_0)}y + B(t_0)) - L(I_0, \rho^{H(t_0)}y + B(t_0))}{\rho^{1-H(t_0)}} \right\|_{L^2(\Omega)} dx. \quad (36)$$

Using essentially the same arguments as those used in proving expression (22) and Remark 3.7 in Boufoussi *et al.* [6], but for the process $X(t) = B(t) - B(t_0)$ instead of the mBm, we obtain

$$\left\| \frac{L(I_0, \theta(\rho)x + \rho^{H(t_0)}y + B(t_0)) - L(I_0, \rho^{H(t_0)}y + B(t_0))}{\rho^{1-H(t_0)}} \right\|_{L^2(\Omega)}$$
$$\leq Ct^{1-H(t_0)(1+\xi)} |x|^\xi \left( \frac{\theta(\rho)}{\rho^{H(t_0)}} \right)^\xi,$$

for sufficently small $\rho$ and all $0 < \xi < \frac{1}{2\sup_{t\geq 0} H(t)} - \frac{1}{2}$. Hence, (36) is dominated by

$$Ct^{1-H(t_0)(1+\xi)} \int_{\mathbb{R}} |f(x)||x|^\xi \, dx \times \left( \frac{\theta(\rho)}{\rho^{H(t_0)}} \right)^\xi. \quad (37)$$

This last integral is finite by assumption (33) and then (37) tends to zero as $\rho$ tends to zero. This completes the proof. □

## Acknowledgments

Research supported by the program Volubilis: "Action intgre MA/06/142."